\documentstyle[12pt]{article}

  \newcommand{\const}{\rm const}

  \newcommand{\argmax}{\rm  argmax}
  \newcommand{\Rad}{\rm Rad}

\begin{document}

 \begin{center}

 \ {\bf Sharp bounds for  the spherical restriction Fourier transform}\par

\vspace{4mm}

\ {\bf in classical  Lebesgue - Riesz and Grand Lebesgue spaces}\par

\vspace{4mm}

{\bf for  ordinary and radial functions.} \par

\vspace{5mm}

 {\bf M.R.Formica,  E.Ostrovsky and L.Sirota.}

 \end{center}

\vspace{4mm}

 \ Universit\`{a} degli Studi di Napoli Parthenope, via Generale Parisi 13, Palazzo Pacanowsky, 80132,
Napoli, Italy. \\

e-mail: mara.formica@uniparthenope.it \\

\vspace{4mm}

 \ Israel,  Bar - Ilan University, department  of Mathematic and Statistics, 59200, \\

\vspace{4mm}

e-mails: \  eugostrovsky@list.ru \\
sirota3@bezeqint.net \\

\vspace{5mm}

 \begin{center}

 \ {\sc Abstract.}

 \end{center}

\vspace{4mm}

  \ We derive bilateral  estimates for the constants appearing in the Fourier transform restricted  theorems on the Euclidean
  sphere for the ordinary and especially {\it radial} functions  belonging to  the Lebesgue - Riesz spaces as well as belonging
  to the Grand Lebesgue Spaces.  \par
  \ We obtain an exact estimate for the norm  of the restriction Fourier transform operator acting on the radial functions. \par

 \vspace{5mm}

  \par

\vspace{5mm}

\section{Definitions. Statement of problem. Previous results.}

\vspace{5mm}

\begin{center}

{\sc Fourier transform  and restriction problem .} \\

\end{center}

 \ Define as ordinary the Fourier transform $ \ \hat{f}(k) = \hat{f}[f](k), \ k \in R^d \ $
 for the integrable function on the whole Euclidean space $ \ R^d, \ d = 2,3,\ldots \ $  as follows

\begin{equation} \label{Fourier}
\hat{f}(k) = \hat{f}[f](k) \stackrel{def}{=} \int_{R^d} e^{-i (k,x)} \ f(x) \ dx.
\end{equation}

 \ Henceforth $ \ (k,x) $ denotes as usually the inner (scalar) product  $ \ k x = \sum_{j=1}^d k_j x_j,  \  r = |x|: =
  \sqrt{x,x}, \ s = |k| := \sqrt{k,k}. \ $ \par

\vspace{4mm}

 \ Denote as ordinary by $ \ ||f||_p = ||f||L_p \ $ the classical Lebesgue - Riesz norm for the  measurable function $ \ f: R^d \to R \ $
 or $ \ f: R^d \to C, \ \ C $ is complex plane,  over the whole space

$$
||f||_p = ||f||L_p(R^d) := \left[ \ \int_{R^d} |f(x)|^p \ dx \ \right]^{1/p}.
$$
 \ Denote also the usually {\it surface }  Euclidean measure defined on the measurable subsets of the unit sphere

$$
  S(d-1)  := \{ x, \ x \in R^d, \ |x| = 1 \},
$$
by $ \ \gamma. \ $ For instance,

$$
A(d) := \gamma(S(d-1)) = \frac{ d \ \pi^{d/2}}{\Gamma(1 + d/2)} = \frac{2 \ \pi^{d/2}}{\Gamma(d/2)}
$$
is the area of sphere $ \ S(d-1). \ $ As usually,  $ \ \Gamma(\cdot) \ $ denotes the Euler's Gamma function. \par
 \ The ordinary Lebesgue - Riesz surface norm
 $ \ ||g||L_q(S(d-1)), \ q \ge 1 \ $ for the  measurable function $ \ g: \ S(d-1) \to C \ $ is defined alike one on the whole space

 \begin{equation} \label{surface norm}
 ||g||L_q(S(d-1)) \stackrel{def}{=} \left[ \ \int_{S(d-1)} |g(y)|^q \ \gamma(dy) \ \right]^{1/q}, \ q \ge 1.
 \end{equation}

\vspace{5mm}

 \hspace{2mm} {\bf Definition 1.1. } The inequality of the form

\begin{equation} \label{key ineq}
||\hat{f}||L_q(S(d-1)) \le K_d(p,q) \ ||f||L_p(R^d),
\end{equation}
with finite coefficient $ \ K_d(p,q) \ $ which holds true for arbitrary $ \ f \in L_p(R^d) \ $ and for some non - trivial domain of values $ \ (p,q) \in B, \ $
is named {\it Fourier restriction inequality.}\par
 \ Here the  non - trivial domain  $ \ B = B(d)  \ $ consists on the whole set of the values of the parameters $ \ (p,q), \ p,q \ge 1,   \  $ for which
$ \ K_d(p,q) < \infty. \ $ \par

 \ Another name: Tomas - Stein inequality, see an pioneer work  \cite{Tomas}, 1975; in this article was obtained  in particular the following important necessary
conditions for  this inequality

\begin{equation} \label{necess cond}
1 \le p \le  \frac{2d +2}{d + 3}, \hspace{3mm} q \le \frac{d-1}{d+1} \cdot \frac{p}{p-1}.
\end{equation}

\vspace{4mm}

 \ There are huge works devoted to the calculation, or at last evaluation, of the domain $ \ B, \ $ as well as the  {\it upper} estimation the coefficient $ \ K_d(p,q); \ $
on the other words, the norm of correspondent restriction  Fourier operator on some regular surfaces, in particular on the sphere,  see e.g. in the works
\cite{Bloom}, \cite{Ferreyra}, \cite{Foschi 1} - \cite{Foschi 2}, \cite{Frank}, \cite{Laura}, \cite{Oliveira 1}  - \cite{Oliveira 2}, \cite{Tomas}, \cite{Tao 1} - \cite{Tao 3},
\cite{Urciuolo}  etc.\par

\vspace{4mm}

 \hspace{3mm} {\bf We intent in this short report to offer the simple lower bounds for this coefficients,   calculate its exact value only for the radial functions,
 find in this case the extremal functions. } \par

\vspace{3mm}

  \ {\bf We extent also these results on the Grand Lebesgue Spaces instead ordinary Lebesgue - Riesz ones.} \par

\vspace{4mm}

 \  We can and will to take as the value $ \ K_d(p,q)  \ $ its minimal value:

\begin{equation} \label{def Z}
K_d(p,q) \stackrel{def}{=} \sup_{ f: ||f||_p = 1}  Z_d(||\hat{f}||L_q(S(d-1) ), ||f||_p),
\end{equation}
where

\begin{equation} \label{minimal value}
Z_d(||\hat{f}||L_q(S(d-1) ), ||f||_p) = Z_d(f) \stackrel{def}{=} \ \frac{||\hat{f}||L_q(S(d-1))}{||f||_p}, \ (p,q) \in B,
\end{equation}

\vspace{3mm}

 and separately only for radial functions

\vspace{3mm}

\begin{equation} \label{def Z Ead}
K_d^{Rad}(p,q) \stackrel{def}{=} \sup_{ f: \ f \in \Rad, \ ||f||_p = 1}  Z_d(||\hat{f}||L_q(S(d-1) ), ||f||_p),
\end{equation}

\vspace{4mm}

 \ Introduce also an extremal function,   which is not uniquely determined, e.g. up to change of sign and up to permutation
 of arguments

\begin{equation} \label{extremal fun}
g(x)  = g_d[p,q](x) \stackrel{def}{=} \argmax_{g: \ ||g||_p \le 1}  \ Z_d(||\hat{g}||L_q(S(d-1) ), ||g||_p),
\end{equation}
so that evidently $ \ ||g||_p =  ||g_d[p,q]||_p  = 1 \ $ and

$$
 Z_d(||\hat{g}||L_q(S(d-1) ), ||g||_p) = K_d(p,q), \ (p,q) \in B,
$$
and analogously for the radial functions:

\begin{equation} \label{extremal rad fun}
g^{\Rad}(x)  = g_d^{\Rad}[p,q](x) \stackrel{def}{=} \argmax_{g: \ g \in \Rad, \ ||g||_p \le 1}  \ Z_d(||\hat{g}||L_q(S(d-1) ), ||g||_p),
\end{equation}
so that evidently $ \ ||g^{\Rad}||_p =  ||g_d^{\Rad}[p,q]||_p  = 1 \ $ and of course

$$
 Z_d^{\Rad}(||\hat{g}^{\Rad}||L_q(S(d-1) ), ||g||^{\Rad}_p) = K_d^{\Rad}(p,q), \ (p,q) \in B.
$$

\vspace{5mm}

\begin{center}

 {\sc About radial functions.   }

\end{center}

\vspace{5mm}

 \ Recall that the function $ \ f = f(x), \ x \in R^d \ $ is said to be {\it radial,} or equally spherical symmetry,
 iff it dependent only  on the polar radii (Euclidean norm) of the argument vector  $ \ x: \ r = |x|: \ $

\begin{equation} \label{radial function}
\exists F = F(r), \ f(x) = F(r) = F(|x|), \ x \in R^d; \ r \in (0,\infty).
\end{equation}

 \ Write: $ \  f(\cdot) \in \Rad. \ $ If this function $ \ f(\cdot) \ $ is radial, then   $ \ \hat{f}[f](k) \ $ is one, as well. \par
 \ More detail, introduce the following kernel function

\begin{equation} \label{kernel fun}
V_d(s,r) \stackrel{def}{=} (2 \pi)^{d/2} \ J_{(d-2)/2}( s \ r) \ s^{(2 - d)/2} \ r^{d/2},
\end{equation}
where  $ \ J_l(\cdot) \ $ denotes as usually the Bessel's function  of order $ \ l; \ l \ge 0. \ $\par
 \ Then the function $ \ \hat{f}(k) \ $ has a form $ \ \hat{f}(k) = G(s) = G(|k|), \ $ where

\begin{equation} \label{Fourier transform}
G(s) = \int_0^{\infty} V_d(s, r) \ F(r) \ dr.
\end{equation}

 \vspace{4mm}

 \ In particular, in this radial case

\begin{equation} \label{full int}
\int_{R^d} f(k) \ dk = \int_{R^d} F(|x|) \ dx = G(0+)  = \frac{2}{\Gamma(d/2)} \int_0^{\infty} r^{d-1} \ F(r) \ dr.
\end{equation}

and

\begin{equation} \label{Lp radial norm}
\int_{R^d} |F(|x|)|^p \ dx = \frac{2}{\Gamma(d/2)} \ \int_0^{\infty} r^{d-1} \ |F(r)|^p \ dr.
\end{equation}

\vspace{5mm}

\section{Main result: lower estimate.}

\vspace{5mm}

 \hspace{3mm}  Let us consider the following example (Gaussian radial density function)

\begin{equation} \label{Gaussian density}
 h(x) = h_{\sigma}(x) \stackrel{def}{=} (2 \pi)^{-d/2} \ \sigma^{-d} \ \exp \left\{ \ - 0.5 \ \sigma^{-2} \ ||x||^2  \ \right\}, \ \sigma = \const \in (0,\infty).
\end{equation}
 \ The Fourier transform of this function has a form

$$
\hat{h}(k) = \exp \left\{ \  - 0.5 \ \sigma^2 \ ||k||^2  \ \right\}, \ k \in R^d.
$$

 \ Note that if $ \  k \in  S(d-1),  \ $ then $ \ \hat{h}(k) = \exp \left\{ \ - 0.5 \ \sigma^2 \ \right\}, \ $ so that this function is constant on the
 surface of the unit sphere $ \ S(d-1) \ $ and following

$$
||\hat{h}||_{q, S(d-1)} = \exp(-\sigma^2/2) \ A^{1/q}(d-1), \ q \ge 1.
$$

\vspace{3mm}

 \ Further,

$$
||h||_p = (2 \pi)^{- 0.5 d (1 - 1/p)} \ \sigma^{-d(1 - 1/p)} \ p^{-d/(2p)}, \ p \ge 1.
$$

\ Therefore, for all the positive values $ \ \sigma \ $

$$
K_d(p,q) \ge e^{-\sigma^2/2} \ A^{1/q}(d-1) \ (2 \pi)^{0.5 d (1 - 1/p)} \ p^{d/(2p)} \ \sigma^{d(1 - 1/p)},
$$
and we conclude after maximization over $ \ \sigma \ $ \par

\vspace{4mm}

{\bf Proposition 2.1.}  We deduce that for all the values of parameters for which

$$
d \ge 2, \ p,q \ge 1 \ \Rightarrow \ K_d(p,q) \ge
$$

 \vspace{3mm}

\begin{equation} \label{lower Gaussian}
 A^{1/q}(d-1) \ (2 \pi)^{0.5 d (1 - 1/p)} \ p^{d/(2p)} \ \times [d(1 - 1/p)]^{0.5 d (1 - 1/p)}.
\end{equation}

 \vspace{5mm}

\section{Main result: sharp  radial estimate.}

\vspace{5mm}

 \hspace{3mm}  Note first of all that there ara some reasons to consider the case of radial functions $ \ f = f(|x|). \ $   Indeed,
 denote by $ \  U \ $ an arbitrary unitary linear operator acting from the space $ \ R^d \ $ into oneself.  It follows immediately
 from the expression for this function  (\ref{extremal fun}) after changes of variables that  $ \ g(Ux) = g(x).  \ $  \par
  \ In detail,  introduce the so - called {\it  unitary dilation operator} $ \ T_U[g](x) \ $ for the function $ \ g: R^d \to R \ $  and for
 arbitrary unitary linear operator  $ \ U: R^d \to R^d \ $

 \begin{equation} \label{unit dilat oper}
 T_U[g](x) \stackrel{def}{=} g(Ux), \ x \in R^d.
 \end{equation}
  \ We have

$$
||T_U[g]||_p^p  = \int_{R^d} |g(Ux)|^p \ dx = \int_{R^d} |g(x)|^p \ dx = ||g||_p^p;
$$
 and quite alike

$$
|| \hat{ T_U g  }||L(q, S(d-1)) = || \hat{ g  }||L(q, S(d-1)).
$$

 \ Following, if the function $ \ g(\cdot), \ ||g||_p \le 1 \ $ is extremal and {\it if she was the only one,}  its unitary  delation $ \ T_U [g] \ $
 is also. Therefore, under this hypothesis (which is false!)

 $$
 T_U[g] = \pm g, \ \Leftrightarrow g(\cdot) \in \Rad.
 $$
 \ We can  and will assume without loss of generality that $ \  T_U[g] = g, \ \Rightarrow g(x) = g_o(|x|) \ $  for some (measurable) numerical valued radial
 function $ \ g_o = g_o(r), \ r > 0. \ $    \par

\vspace{3mm}

 \ Thus, the function  $ \ g(x) \ $ is radial, and one can find the extremal in  (\ref{def Z}), (\ref{minimal value}) only among the set of
 all radial functions: $ \ g(x)  = F(|x|) = F(r), \ r > 0.  \ $  \par

 \ Note that here

 \begin{equation} \label{p root spherical}
||F||_p \  = \sqrt[p]{\frac{2}{\Gamma(d/2)} \ \left[ \ \int_0^{\infty} r^{d - 1} \ |F(r)|^p \ dr  \right] }, \  p > 0.
 \end{equation}

 \ Further, it follows from (\ref{Fourier transform})

\begin{equation} \label{Fourier  transform  unit sphere}
G(s) = \int_0^{\infty} V_d(r) \ F(r) \ dr, \ |s| = 1,
\end{equation}
i.e. this function is constant on the surface of unit sphere $ \ S(d-1). \ $ Obviously,

$$
V_d(r) = (2 \pi)^{d/2} \ J_{(d - 2)/2}(r) \ r^{d/2}, \ r > 0.
$$
 \ Therefore

$$
|| G(\cdot) ||_{q, S(d-1)} = A^{1/q}(d-1) \int_0^{\infty} V_d(r) \ F(r) \ dr.
$$

\vspace{5mm}

 \ Thus, the investigated here problem is reduced to the following one: \par

\vspace{4mm}

\begin{equation} \label{def Z0}
W_{d,p,q}[F] \stackrel{def}{=} \left[ \ 2^{-1/p} \ \Gamma^{1/p}(d/2) \ A^{1/q}(d-1) \ \right] \times
\end{equation}

\begin{equation} \label{def ZZ}
\frac{\int_0^{\infty} V_d(r) \ F(r) \ dr }{ \left[ \ \int_0^{\infty} r^{d-1} |F(r)|^p \ dr  \ \right]^{1/p}} \to \sup / \{F \in L(p, R_+, r^{d-1} dr) \}.
\end{equation}

  \ This problem may be transformed as follows. Introduce the following measure $ \ d \nu = \nu_d(dr) := r^{d-1} \ dr;   \  $ then we get to
the following problem

\begin{equation} \label{key problem}
\Phi[F] = \Phi_{d,p,q}[F] \stackrel{def}{=} \frac{\int_0^{\infty} r^{1 - d} \ V_d(r) \ F(r) \ \nu_d(dr)}{||F||_{p, \nu}}
 \to \sup / \{F \in L(p,\nu) \}.
\end{equation}

 \ Denote as usually $ \ p' = p/(p-1), \ p \in (1,\infty); \   1' := +\infty;   \ \infty' := 1.\ $  If follows immediately from  the H\"older's
inequality that

\begin{equation} \label{main result}
Q[d,p,q] \stackrel{def}{=} \sup_{F \in L(p)} \Phi_{d,p,q}[F] = \left\{ \ \int_0^{\infty} r^{p' (1 - d)} \ V_d^{p'}(r) \ \nu_d(dr) \ \right\}^{1/p'}.
\end{equation}

\vspace{4mm}

 \ Let us introduce the following coefficient

\begin{equation} \label{coeff P}
P = P(d; p,q) \stackrel{def}{=} 2^{-1/p} \ \Gamma^{1/p}(d/2) \ A^{1/q}(d-1) \  (2\pi)^{-d/(2 p')}.
\end{equation}

  \ To summarize.\par

  \vspace{4mm}

 \ {\bf Theorem 3.1.}

\begin{equation} \label{rad value}
K^{\Rad}(d; p,q) = P(d;p,q) \cdot \left\{ \ \int_0^{\infty} r^{ (2 + d(p-2))/( 2(p-1) ) } \ J^{p'}_{(d-2)/2}(r) \ dr \  \right\}^{1/p'}.
\end{equation}

\ Herewith the  arbitrary  extremal function $ \ F_o(r) \ $ has a form

\begin{equation} \label{extremal fun}
F_0(r) = C \cdot [V_d(r)]^{1/(p-1)}, \ p > 1, \ C = \const.
\end{equation}

\vspace{4mm}

 \ Notice that the integral in (\ref{rad value}) convergent under our  assumptions iff

\begin{equation} \label{condit converg}
1 < p < \frac{2d}{d + 1}.
\end{equation}

\ As regards for the properties of Bessel's functions see, e.g. \cite{Watson}. Namely, as $ \ z \to 0+ \ $

$$
J_{\alpha}(z) \asymp C_0(\alpha) \ z^{\alpha},
$$
 and correspondingly as $ \ z \to \infty \ $

$$
J_{\alpha}(z) \asymp   C_{\infty}(\alpha) \ z^{-1/2}, \ \alpha = \const \ge 0.
$$

\vspace{5mm}

\section{Generalization on the Grand Lebesgue Spaces.}

\vspace{5mm}

 \hspace{3mm} Let $ \ (a,b) = \const, \ 1 \le a < b \le \infty, \ $ and let $ \ \psi = \psi(p), \ p \in (a,b) \ $  be
 bounded from below: $ \ \inf_{p \in (a,b) } \psi(p) > 0 \ $ measurable function. The set of all such a functions will be
 denoted by  $ \ \Psi(a,b); \ $ put also

$$
\Psi  := \cup_{(a,b): 1 < a  < b < \infty} \Psi(a,b).
$$

\vspace{4mm}

 \ {\bf Definition  4.1.} The Grand Lebesgue Space $ \ G\psi, \ \psi \in \Psi(a,b) \ $  builded over the set $ \ R^d, \ $ or equally over the sphere
 $ \ S(d-1), \ $ consists by definition on all the integrable functions having  a finite norm

\begin{equation} \label{GLS norm Rd}
||f||G\psi(R^d) \stackrel{def}{=} \sup_{p \in (a,b)} \left\{ \ \frac{||f||_p}{\psi(p)}   \ \right\}.
\end{equation}

\vspace{4mm}

 \ These space was investigated in many works, see e.g. \cite{Buldygin}, \cite{Ermakov etc. 1986},
\cite{Fiorenza 1}, \cite{Fiorenza 2}, \cite{Fiorenza-Formica-Gogatishvili-DEA2018},
\cite{fioforgogakoparakoNAtoappear}, \cite{fioformicarakodie2017}, \cite{formicagiovamjom2015}, \cite{Iwaniec},
\cite{Kozachenko 1}, \cite{Kozachenko 2}, \cite{Ostrovsky 0} - \cite{Ostrovsky 3}. In particular, the belonging of
the function to certain Grand Lebesgue Space $ \ G\psi \ $ is closely related with its tail behavior and  is related
with its moment generating function $ \ \int \exp(\lambda f(x)) \ dx. \ $ \par
 \ Suppose now that the function $ \ f = f(x) \ $ belongs to sone GLS $ \ G\psi(a,b), \ \exists \psi \in \Psi: \ $

$$
||f||G\psi \le \psi(p), \ p \in (a,b).
$$
 \ Define the (cut) set

$$
D = \{ \ p, \ \exists q \ \Rightarrow (p,q) \in B \ \},
$$
and assume its non - triviality: $ \ D  \ne \emptyset. \ $ We have

$$
||\hat{f}||_{q,S(d-1)} \le ||f||G\psi \cdot \psi(p) \cdot K(p,q), \ p \in (a,b), \ q \in D.
$$
 \ Therefore

$$
||\hat{f}||_{q,S(d-1)} \le ||f||G\psi \cdot \psi(p) \cdot \zeta(q), \ q \in D,
$$
where

$$
\zeta(q) := \inf_{p \in D} [ \ \psi(p) \cdot K(p,q)  \  ].
$$

 \vspace{3mm}

 \ To summarize: \par

\vspace{4mm}

 \ {\bf Theorem 4.1.}  We deduce  under formulated restrictions

$$
||\hat{f}||_{G\zeta,S(d-1)} \le 1 \times ||f||G\psi,
$$
where the constant $ \ "1" \ $ is the best possible. \par
 \ The non - improvability for this constant is grounded in particular in  \cite{Ostrovsky 1}.

 \vspace{6mm}

\vspace{0.5cm} \emph{Acknowledgement.} {\footnotesize The first
author has been partially supported by the Gruppo Nazionale per
l'Analisi Matematica, la Probabilit\`a e le loro Applicazioni
(GNAMPA) of the Istituto Nazionale di Alta Matematica (INdAM) and by
Universit\`a degli Studi di Napoli Parthenope through the project
\lq\lq sostegno alla Ricerca individuale\rq\rq .\par

\end{document}